\documentclass[a4paper,11pt]{article}
\usepackage{cite}
\usepackage{graphicx}
\usepackage[tight,footnotesize]{subfigure}
\usepackage{amsfonts,latexsym,amssymb,amsmath}
\usepackage{epstopdf}
\usepackage{tikz}
\usetikzlibrary{arrows}
\usetikzlibrary{arrows,patterns,calc,patterns,decorations.pathmorphing,decorations.markings}
\newcommand{\vect}[1]{\mathbf{#1}}

\newcommand{\factor}{q_i\left(s,\,k_1,\,k_2,\,\tau,\,\lambda_i\right)}
\newcommand{\qed}{\hfill\blacksquare}
\newtheorem{lemma}{\bfseries{Lemma}}
\newtheorem{definition}{\bfseries{Definition}}

\begin{document}
\title{The Most Exigent Eigenvalue: Guaranteeing Consensus under an Unknown Communication Topology and Time Delays}
\author{Rudy Cepeda-Gomez, Torsten Jeinsch, and Wolfgang Drewelow\\Institute of Automation, University of Rostock, Rostock, Germany. \\ \textit{rudy.cepeda-gomez@uni-rostock.de,}\\ \textit{torsten.jeinsch@uni-rostock.de, wolfgang.drewelow@uni-rostock.de}}%
\maketitle
\begin{abstract}
In the analysis of consensus problems for multi-agent systems affected by time delays, the delay margin for a given protocol is heavily dependent on the communication topology being used. This document aims to answer the question of what is the largest delay value under which convergence to consensus is still possible even if the communication topology is not known. To answer this question we revisit the concept of most exigent eigenvalue, applying it to two different consensus protocols for agents driven by second order dynamics. We show how the delay margin depends on the structure of the consensus protocol and the communication topology, and arrive to a boundary that guarantees consensus for any connected communication topology. The switching topologies case is also studied and It is shown that for one protocol the stability of the individual topologies is sufficient to guarantee consensus in the switching case, whereas for the other one it is not.
\end{abstract}

\section{Introduction}\label{sec:int}
The coordinated execution of a task by a team of multiple dynamic agents almost always requires that a decision upon the value of a variable is reached. If we talk about formation flight of multiple unmanned aerial vehicles, for example, a common heading and a common speed must be selected. When this common value is to be reached by local interactions among the agents of the group a \emph{consensus} problem appears. The work Vicsek \cite{vicsek1995} was one of the earliest studies of the properties of such agreement problem, considering discrete agents governed by first order dynamics. This work was later expanded Jadbabaie \cite{Jadbabaie2003} and others. It was the work of Olfati-Saber and Murray \cite{Olfati-Saber2004} the one to introduce the term \emph{consensus}. They exhaustively studied the conditions under which a control protocol leads a group of agents governed by first order continuous dynamics to reach an agreement, including cases with switching topologies and time delays. Many derivatives of this work have appeared in the years since its publication, including the extension to second order dynamics \cite{Ren2008,Lin2009-1} with the inclusion of switching topologies \cite{Ren2005}, time delays \cite{Lin2008-1,Lin2008-2,Sun2009-1,Sun2009-2,Muenz2010,Liu2009,Liu2010,Hu2010}, or both \cite{Lin2009-2,Lin2010,Xiao2006,Xiao2008}. A good overview of the current challenges and results in this topic can be found in \cite{Cao2013b,Davoodi2016}.

One of the most important contributions of \cite{Olfati-Saber2004} was  to pinpoint the importance that some results from algebraic graph theory have in the topic of consensus. They showed that for the agents to reach an agreement upon the average of the initial conditions, the graph that describes the communication topology must be connected and balanced. Another important result is that the time constant of the system, i.e., the time the agent need to reach consensus, is related to the second smallest eigenvalue of the Laplacian matrix of the graph, also known as the algebraic connectivity or the Fiedler eigenvalue \cite{Fiedler1973}. They also show that when a time delay is present, the largest eigenvalue of the Laplacian limits the tolerance of the system to a delay. These results led to some works in which the weights of the communication channels are designed such that a certain delay margin is attained \cite{Sipahi2011,Qiao2013,Qiao2014}.

Almost all works dealing with time delayed consensus problems, whether they consider first or second order dynamic models, use control laws in which the delay affects the information of the own state of each agent as well as the state information coming from other agents. This so called \emph{self delay} case can be seen, for example, in the works \cite{Olfati-Saber2004,Lin2008-1,Lin2008-2,Sun2009-1,Sun2009-2,Muenz2010,Liu2009,Liu2010,Hu2010}. The consensus protocols which share this structure are always described using the Laplacian matrix. A different case in which the state of the agent is not affected by the delay, and therefore the Laplacian matrix plays no role is presented in \cite{IJC1}. 

It was shown in \cite{TAC}, however, that consensus protocols with or without self delay have a common structure which can be used to simplify its analysis with respect to the time delay. This paper shows how the characteristic equation of a consensus problem can be decomposed in factors of low order (equal to the order of the dynamics of the individual agents) with a similar structure. The factors differ only in a certain coefficient given by an eigenvalue of a matrix related to the communication structure. In \cite{IJC1}, the authors take advantage of this property to perform an exhaustive analysis of a protocol without self delay for a fixed communication topology. In that paper the concept of \emph{most exigent eigenvalue}, defined as the eigenvalue of the matrix which creates the most restrictive stability boundary with respect to the time delay, is introduced.

In the present paper, we aim to perform a similar analysis for a consensus protocol with self delay. We show how the most exigent eigenvalue in that case is the largest eigenvalue of the Laplacian matrix. We compare this to our previous result for the case without self delay and show how it can be used to guarantee consensus under an unknown communication topology in the presence of time delays. Furthermore, we study the switching topologies case and show how the protocol with self delay is stable regardless of the switching scheme, whereas the protocol with self delay can be destabilized if the communication topologies switches in the correct pattern.

Although the term \emph{most exigent eigenvalue} was first introduced by \cite{IJC1}, the work \cite{Olfati-Saber2004} already showed how there is one eigenvalue setting the most restrictive delay boundary for a given communication topology, and \cite{Sipahi2011} extended this. These works, however, are focused on the case of agents with first order dynamics, fact that separates them from the present article, in which second order agents are considered.

It is important to highlight that, when it comes to the stability analysis of consensus systems with respect to time delays, almost all the previous works rely on Lyapunov–Krasovskii- or Razhumikin-based methodologies \cite{Lin2008-1,Lin2008-2,Lin2009-1,Lin2009-2,Meng2011,Zhang2014} or a generalized Nyquist criterion \cite{Liu2009,Muenz2010}. All of these treatments provide only sufficient conditions on the delays to achieve asymptotic stability. They only produce very conservative results, leading to stability bounds for very small delays.

As a different way to analyze the consensus problem with respect to the delay,  a methodology for the analysis of consensus protocols with single and multiple, rationally independent time delays was introduced in \cite{TAC,IJC2}. This earlier work is based on the combination of a factorization procedure and the deployment of the Cluster Treatment of Characteristic Roots methodology \cite{Ergenc2007,Fazelinia2007,Olgac2005}, which provides exact, exhaustive and explicit stability regions in the domain of the delays. 

The paper is organized as follows. Section \ref{sec:prot} presents the two consensus protocols under study and some brief results regarding their stability analysis. Section \ref{sec:mee} presents the results on the concept of most exigent eigenvalue for both protocols and Section \ref{sec:switch} studies the switching topologies case. Some conclusions and directions for future work are presented in Section \ref{sec:conc}. In the rest of the paper, scalars are represented by lowercase italic letters ($k,\,\lambda,\,\tau$), vectors by lowercase bold letters ($\mathbf{x}$) and matrices by uppercase boldface letters ($\mathbf{A},\boldsymbol{\Lambda}$).

\section{Consensus Protocols and Stability Analysis}\label{sec:prot}
In this work we consider a set of $n$ one dimensional agents driven by second order dynamics\footnote{This double integrator model is considered instead of a more general dynamic model to simplify the stability analysis with respect to the delay. An example of such stability analysis for a more general dynamic system can be found in \cite{IJSS2}.} of the form $\ddot{x}_i(t)=u_i(t)$, $i=1,\,2\,\ldots,\,n$. The control input $u_i(t)$ is computed based on the state of agent number $i$ as well as the state of some of its peers, known as the informers of $i$. The set of informers of agent $i$ is denoted as $\mathcal{N}_i$, and its cardinality, \emph{i.e.}, the number of informers of agent $i$, is denoted as $\delta_i$. For this study we consider that the communication is bi-directional, \emph{i.e.}, $i\in\mathcal{N}_j \Leftrightarrow j\in\mathcal{N}_i$. This means that the network can be represented by an undirected graph. The \emph{adjacency matrix} of this graph, denoted as $\vect{A}_\Gamma=[a_{ij}]$, is defined such that $a_{ik}=a_{ki}=1$ when agents $i$ and $j$ share a communication link and $a_{ik}=a_{ki}=0$ otherwise. The diagonal elements of $\mathbf{A}_\Gamma$ are considered to be zero. The \emph{degree matrix} of the communication topology is a diagonal matrix $\boldsymbol{\Delta}$ with its $i$-th diagonal element equal to the degree of agent $i$. The \emph{Laplacian matrix} of the communication topology is defined as $\vect{L}=\boldsymbol{\Delta}-\mathbf{A}$, and the \emph{weighted adjacency matrix} is $\vect{C}=\boldsymbol{\Delta}^{-1}\mathbf{A}_\Gamma$. We also assume that the communication among agents is corrupted by a time delay $\tau$, which is constant and uniform across the network.

The uniformity assumption for the communication delay is a common one in the field of consensus and is based on the fact that members of the multi-agent systems under consideration are identical and have therefore similarly limited bandwidths in their sensing or communication capabilities. Cases in which a subset of the agents or even each individual agent had their own time delay could still be considered, but some mathematical simplifications could not be used in that case, making the stability analysis prohibitively complex.

In this work we study two different control actions that can be used by the agents in order to reach consensus. The first one was introduced in \cite{TAC,IJC1} and is given as
\begin{equation}
\begin{split}
u_i\left(t\right)=&k_1\left[\sum_{j\in\mathcal{N}_i}{\left(\frac{x_j\left(t-\tau\right)}{\delta_i}\right)}-x_i\left(t\right)\right]\\+&k_2\left[\sum_{j\in\mathcal{N}_i}{\left(\frac{\dot{x}_j\left(t-\tau\right)}{\delta_i}\right)}-\dot{x}_i\left(t\right)\right],
\end{split}
\label{eq:noself}
\end{equation}
whereas the second one, studied in \cite{Lin2008-1}, is
\begin{equation}
\begin{split}
u_i\left(t\right)=&k_1\left[\sum_{j\in\mathcal{N}_i}{\left(x_j\left(t-\tau\right)-x_i\left(t-\tau\right)\right)}\right]\\+&k_2\left[\sum_{j\in\mathcal{N}_i}{\left(\dot{x}_j\left(t-\tau\right)-\dot{x}_i\left(t-\tau\right)\right)}\right].
\end{split}
\label{eq:self}
\end{equation}
In equations \eqref{eq:noself} and \eqref{eq:self} $k_1$ and $k_2$ are positive control gains selected by the user. For $\tau=0$ both protocols guarantee consensus provided that the communication topology is connected. Connectivity is therefore assumed as granted for the rest of the paper.

The main difference between protocols \eqref{eq:noself} and \eqref{eq:self} lies in the presence of the so called \emph{self delay}. In \eqref{eq:noself}, agent $i$ uses the delayed information coming from its informers and compares it to its own current state. In \eqref{eq:self} the delay is present in the state of agent $i$ as well as in the state of the informers. This structural difference has some implications in the dynamics of the group, a fact that is presented in the following paragraphs.

Despite the structural differences, it was shown in \cite{TAC} that the characteristic equations of both protocols can be expressed as a product of $n$ second order factors. This factorization property stems from the fact that both protocols can be represented in state space as
\begin{equation}
\dot{\vect{x}}\left(t\right)=\left(\mathbf{I}_n\otimes\vect{F}_1\right)\vect{x}\left(t\right)+\left(\mathbf{M}\otimes\vect{F}_2\right)\vect{x}\left(t-\tau\right),
\label{eq:genss}
\end{equation}  
with $\vect{x}=\left[x_1\,\dot{x}_1\,x_2\,\dot{x}_2\,\cdots x_n\,\dot{x}_n\right]^T$, $\mathbf{I}_n$ being the identity matrix of order $n$ and $\otimes$ the Kronecker product. The matrices $\vect{F}_{1,2}$ depend on the particular protocol and the control gains, whereas the matrix $\vect{M}$ depends on the protocol and the communication topology. For \eqref{eq:noself} we have
\begin{equation}
\begin{split}
&\vect{F}_1=\left[\begin{array}{rr}0&1\\-k_1&-k_2\end{array}\right]\ \vect{F}_2=\left[\begin{array}{rr}0&0\\k_1&k_2\end{array}\right]\\ &\vect{M}=\boldsymbol{\Delta}^{-1}\vect{A}_{\Gamma}=\vect{C},
\end{split}
\end{equation}
whereas for \eqref{eq:self}
\begin{equation}
\begin{split}
&\vect{F}_1=\left[\begin{array}{rr}0&1\\0&0\end{array}\right]\ \vect{F}_2=\left[\begin{array}{rr}0&0\\k_1&k_2\end{array}\right]\\ &\vect{M}=\vect{A}_{\Gamma}-\boldsymbol{\Delta}=-\vect{L}.
\end{split}
\end{equation}

Since we are considering undirected graphs, the matrix $\mathbf{M}$ of each protocol is always diagonalizable and its eigenvalues are always real. This is because the weighted adjacency matrix $\mathbf{C}$ is a symmetrizable matrix \cite{Sergienko2003} and the Laplacian $\mathbf{L}$ is symmetric. As proposed in \cite{TAC}, we take advantage of this fact to define the state transformation $\boldsymbol{\xi}=\left(\vect{T}^{-1}\otimes\vect{I}_2\right)\vect{x}$, where $\vect{T}$ is the matrix that diagonalizes $\vect{M}$, \emph{i.e}, $\vect{T}^{-1}\vect{MT}=\boldsymbol{\Lambda}$ with $\boldsymbol{\Lambda}$ being diagonal. Under this transformation, the dynamics of the system in \eqref{eq:genss} are expressed as
\begin{equation}
\dot{\boldsymbol{\xi}}\left(t\right)=\left(\mathbf{I}_n\otimes\vect{F}_1\right)\boldsymbol{\xi}\left(t\right)+\left(\boldsymbol{\Lambda}\otimes\vect{F}_2\right)\boldsymbol{\xi}\left(t-\tau\right).
\label{eq:transformed}
\end{equation}
Given that $\boldsymbol{\Lambda}$ is a diagonal matrix, the system \eqref{eq:transformed} can be seen as a set of $n$ decoupled second order systems, each one of the form
\begin{equation}
\left[\begin{array}{c}\dot{\xi}_i(t)\\\ddot{\xi}_i(t)\end{array}\right]=\vect{F}_1\left[\begin{array}{c}\xi_i(t)\\\dot{\xi}_i(t)\end{array}\right]+\lambda_i\vect{F}_2\left[\begin{array}{c}\xi_i(t-\tau)\\\dot{\xi}_i(t-\tau)\end{array}\right]\quad i=1,\,2,\,\ldots,\,n,
\label{eq:subsys}
\end{equation}
in which $\lambda_i$, $i=1,\,2\,\ldots,\,n$, are the eigenvalues of the matrix $\mathbf{M}$ corresponding to each protocol, which, as stated before, are all real. Furthermore, since $\mathbf{C}$ is a hollow stochastic matrix \cite{Marcus1996}, Gershgorin's circle theorem \cite{Bell1965} guarantees that $\left|\lambda_i\right|\leq 1$ for protocol \eqref{eq:noself}; whereas $\lambda_i\leq0$ for protocol \eqref{eq:self}, because the eigenvalues of the Laplacian are always nonnegative \cite{Biggs1993}.

The characteristic equation of the whole system can be expressed as the product of the characteristic equation of the individual subsystems in \eqref{eq:transformed}. That is, the stability of the system is determined by the roots an equation expressed as the product of $n$ second order factors
\begin{equation}
\begin{split}
CE\left(s,\tau\right)=&\prod_{i=1}^n\factor=\\&\prod_{i=1}^n{\det\left(s\vect{I}_2-\vect{F}_1-\lambda_i \vect{F}_2e^{-\tau\,s}\right)}.
\end{split}
\label{eq:ce}
\end{equation}
Each one of the subsystems in \eqref{eq:subsys} or, equivalently, each one of the factors in \eqref{eq:ce}, corresponds to the dynamics of a certain linear combination of the positions of the agents in the group. These are the states denoted by $\xi_i(t)$.

\emph{Remark:} Complex conjugate eigenvalues could appear if the communication between agents is not bi-directional, as assumed here. In this case the, factorization property has to be extended to include factors of order higher thank the dynamics of the individual agents, as presented in \cite{IJC3}.

Since we are assuming that the topology is connected, a special eigenvalue is present in each protocol. This eigenvalue corresponds to an eigenvector of $\mathbf{M}$ in which all the elements are equal. Without loss of generality, we denote this eigenvalue as $\lambda_1$. For \eqref{eq:noself} the special eigenvalue is $\lambda_1=1$, whereas for protocol \eqref{eq:self} it is $\lambda_1=0$. The factor of \eqref{eq:ce} generated by $\lambda_1$ dictates the dynamics of the group decision value. That is, if all the other factors are stable the agents reach consensus and move together in a trajectory dictated by the dynamics of this factor. We call it the \emph{centroid factor}. The other $n-1$ factors define whether the agents reach an agreement or not, and they are called the \emph{disagreement factors}. 

It was shown in \cite{IJC1} that for protocol \eqref{eq:noself} the group decision value is dictated by a weighted average, defined as
\begin{equation}
\xi_1\left(t\right)=\frac{\sum_{i=1}^n{\delta_ix_i(t)}}{\sum_{i=1}^n{\delta_i}},
\label{eq:weightedc}
\end{equation} 
for which the dynamics is given by
\begin{equation}
s^2+\left(k_2s+k_1\right)\left(1-e^{\tau\,s}\right)=0.
\label{eq:centdyn}
\end{equation}
Notice that \eqref{eq:centdyn} has a root at $s=0$ for any value of $\tau$. This means that this factor is at best marginally stable, which implies that consensus, if reached, will be at a constant position with zero velocity.

For protocol \eqref{eq:self}, on the other hand, previous works \cite{Lin2008-1,CDC11} showed that the agents reach an average consensus given by
\begin{equation}
\xi_1\left(t\right)=\frac{\sum_{i=1}^n{x_i(t)}}{n},
\label{eq:c}
\end{equation} 
with dynamics
\begin{equation}
s^2=0.
\label{eq:avedyn}
\end{equation}
Equation \eqref{eq:avedyn} has two roots at the origin, which do not depend on the delay. This indicates that protocol \eqref{eq:self}, when it is stable, guides the agents towards a constant velocity with linearly increasing position.

For the agents to reach consensus using protocol \eqref{eq:noself} $n-1$ disagreement factors of the form
\begin{equation}
\factor=s^2+\left(k_2s+k_1\right)\left(1-\lambda_ie^{-\tau\,s}\right)=0,
\label{eq:factor1}
\end{equation}
with $i=2,\,3,\,\ldots,\,n$, must be stable. When protocol \eqref{eq:self} is used, the $n-1$ disagreement factors have the form\footnote{Notice that a change in the sign of $\lambda_i$ has been introduced in \eqref{eq:factor2}. This is to use the eigenvalues of $\vect{L}$, which are all positive, and not those of $-\vect{L}$.}:
\begin{equation}
\factor=s^2+\lambda_i\left(k_2s-k_1\right)e^{-\tau\,s}=0.
\label{eq:factor2}
\end{equation}

Factors \eqref{eq:factor1} and \eqref{eq:factor2} are stable for $\tau=0$. As $\tau$ increases, the characteristic roots of the factors move in the complex plane until the delay reaches a critical value for which a root crosses to the right of the complex plane, introducing instability. The \emph{Cluster Treatment of Characteristic Roots Paradigm} (CTCR), introduced in \cite{Olgac2002}, is used to detect these destabilizing crossings. In order to limit the length of this paper, the details of the deployment of CTCR to the factors of the form \eqref{eq:factor1} and \eqref{eq:factor2} are left out. They were already presented in \cite{IJC1}  and \cite{CDC11}. We only state here the final results. 

For factors of the form \eqref{eq:factor1}, the first destabilizing root crosses at a frequency
\begin{equation}
\omega^2=k_1-\frac{\mu}{2}k_2^2+\frac{\sqrt{\rho}}{2},
\label{eq:omega1}
\end{equation}
with
\begin{equation}
\mu=1-\lambda_i^2,\quad \text{and}\quad \rho=\left(\mu k_2^2-2k_1\right)^2-4k_1^2\mu,
\end{equation}
and this crossing occurs at a delay value which depends on the sign of $\lambda_i$
\begin{equation}
\begin{aligned}
\tau&=\frac{1}{\omega}\left(\arctan\left(\frac{-k_2\omega^3}{k_1^2+\omega^2\left(D^2-P\right)}\right)\right)\quad\text{ for }\lambda_i>0,\\
\tau&=\frac{1}{\omega}\left(\arctan\left(\frac{-k_2\omega^3}{k_1^2+\omega^2\left(D^2-P\right)}\right)+\pi\right)\quad\text{ for }\lambda_i<0.
\end{aligned}
\label{eq:tau1}
\end{equation}
For protocol \eqref{eq:self}, the expressions for the crossing frequency and delay are
\begin{align}
\omega^2&=\frac{k_2^2\lambda_i^2+\sqrt{k_2^4\lambda_i^4+4k_1^2\lambda_i^2}}{2},\label{eq:omega2}\\
\tau&=\frac{1}{\omega}\arctan\left(\frac{k_2\omega}{k_1}\right).\label{eq:tau2}
\end{align}

The stability analysis technique presented in \cite{TAC} departs from the knowledge of the communication topology. From here, the eigenvalues of the weighted adjacency or the the Laplacian matrix are found and the stability region in the parametric domain $\left(k_1,\,k_2,\,\tau\right)$ is obtained for each factor. These regions are then intersected to find the combinations of parameters which make the complete system stable.

In the next section, we show that an \emph{a priori} knowledge of the communication topology is not needed to find a delay bound.
\section{The most Exigent Eigenvalue}\label{sec:mee}
Section \ref{sec:prot} showed that each eigenvalue of $\vect{M}$ introduces a stability boundary. However, there is always one eigenvalue that introduces the most restrictive of these boundaries with respect to the time delay, defining therefore the global stability region. This eigenvalue is declared the \emph{most exigent eigenvalue}. This section formally defines this concept, first introduced in \cite{IJC1} for a second order protocol, and declares which eigenvalue is the most exigent for protocols \eqref{eq:noself} and \eqref{eq:self}.

\begin{definition}
For a group of agents interacting under any of the consensus protocols defined in \eqref{eq:noself} or \eqref{eq:self}, the \emph{most exigent eigenvalue} is the eigenvalue of the corresponding $\vect{M}$ matrix that generates the particular factor $\factor$ in the characteristic equation \eqref{eq:ce} which introduces the smallest destabilizing crossing as $\tau$ increases starting from 0, for a fixed set of $k_1$ and $k_2$ values.
\end{definition} 

\emph{Remark:} When unidirectional communications are considered and complex conjugate eigenvalues appear, the factors of order higher than the dynamics of the individual agents may present more than one stability crossing for a single parametric combination. A \emph{most exigent eigenvalue} is hard to define in this case, and the discussions are therefore limited to the bi-directional communications case assumed here.

\begin{lemma}\label{le:me1}
For the consensus protocol without self-delay defined by \eqref{eq:noself}, the most exigent eigenvalue is the smallest eigenvalue of the weighted adjacency matrix $\vect{C}$.
\end{lemma}
\emph{Proof:} The proof was presented in \cite{IJC1}.
$\qed$
\begin{lemma}\label{le:me2}
For the consensus protocol with self-delay defined by \eqref{eq:self}, the most exigent eigenvalue is the largest eigenvalue of the Laplacian matrix $\vect{L}$.
\end{lemma}
\emph{Proof:} We prove this following the same path used in \cite{IJC1} to prove Lemma \ref{le:me1}. We show that the following statements are true when we consider $\omega>0$:
\begin{itemize}
\item[(a)] $\omega$ is a monotonic one to one function of $\lambda$ and $d\omega/d\lambda>0$.
\item[(b)] $\tau$ is a monotonic function of $\omega$ and $d\tau/d\omega>0$. 
\item[(c)] $\tau$ is a monotonic function of $\lambda$ and $d\tau/d\lambda<0$.
\end{itemize}

In order to prove (a), we use $\omega^2=\gamma$ and consider \eqref{eq:omega2} for two different $\lambda$, recasting the equations as
\begin{subequations}\label{eq:gamma}
\begin{align}
\gamma^2-\left(k_2\lambda_1\right)^2\gamma-\left(k_1\lambda_1\right)^2&=0,\label{eq:gammaA}\\
\gamma^2-\left(k_2\lambda_2\right)^2\gamma-\left(k_1\lambda_2\right)^2&=0.\label{eq:gammaB}
\end{align}
\end{subequations}
If this two polynomial equations have a common root, the Sylvester resultant matrix
\begin{equation}
\vect{R}=\left[\begin{array}{cccc}
1&-\left(k_2\lambda_1\right)^2&-\left(k_1\lambda_1\right)^2&0\\
0&1&-\left(k_2\lambda_1\right)^2&-\left(k_1\lambda_1\right)^2\\
1&-\left(k_2\lambda_2\right)^2&-\left(k_1\lambda_2\right)^2&0\\
0&1&-\left(k_2\lambda_2\right)^2&-\left(k_1\lambda_2\right)^2
\end{array}\right]
\end{equation}
must be singular. We have that $\det\left(\vect{M}\right)=k_2^4(\lambda_1^2 - \lambda_2^2)^2$, which can be zero if and only if $\lambda_1=\lambda_2$. Therefore, $\gamma$ is a one to one function of $\lambda$. Since $\omega=\sqrt{\gamma}$, the first part of (a) is proven.

By using implicit differentiation, we see from \eqref{eq:gammaA} that
\begin{equation}
\frac{d\gamma}{d\lambda}=\frac{\left(k_2\lambda\right)^2}{2\gamma},
\label{eq:dgamma}
\end{equation}
and combining \eqref{eq:dgamma} with $d\omega/d\gamma=1/(2\sqrt{\gamma})$ we have
\begin{equation}
\frac{d\omega}{d\lambda}=\frac{d\omega}{d\gamma}\frac{d\gamma}{d\lambda}=\frac{\left(k_2\lambda\right)^2}{4\gamma^{3/2}}=\frac{\left(k_2\lambda\right)^2}{4\omega^{3}}>0,
\end{equation}
which proves the second part of (a).

\begin{figure}
\centering
\includegraphics[scale=0.8]{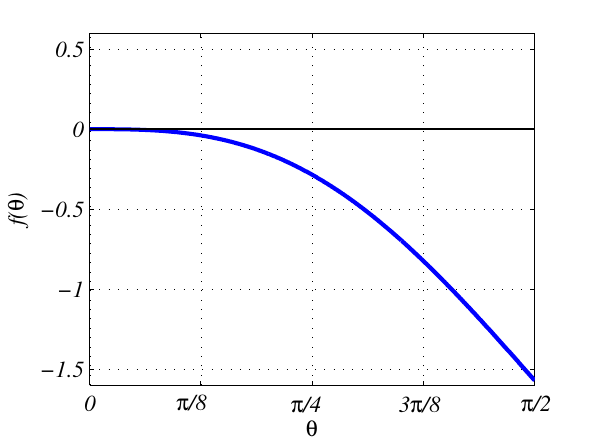}
\caption{Plot of $f(\theta)=\frac{1}{2}\sin(2\theta)-\theta$ for $0\leq\theta\leq\pi/2$.}
\label{fig:1}
\end{figure}

To prove part (b), we consider \eqref{eq:tau2} and take its derivative with respect to $\omega$
\begin{equation}
\frac{d\tau}{d\omega}=\frac{k_2/k_1}{\omega\left(1+\left(\frac{k_2\omega}{k_1}\right)^2\right)}-\frac{1}{\omega^2}\arctan\left(\frac{k_2\omega}{k_1}\right).
\label{eq:dt1}
\end{equation}
We are interested in showing that $d\tau/d\omega<0$. By multiplying \eqref{eq:dt1} times $\omega^2$ we obtain the following inequality
\begin{equation}
\frac{\frac{k_2\omega}{k_1}}{1+\left(\frac{k_2\omega}{k_1}\right)^2}-\arctan\left(\frac{k_2\omega}{k_1}\right)<0.
\label{eq:dt2}
\end{equation}
Using the substitution $\frac{k_2\omega}{k_1}=\tan(\theta)$, valid for $0<\theta<\pi/2$, we transform \eqref{eq:dt2} into
\begin{equation}
\frac{\tan\theta}{1+\tan^2\theta}-\theta<0,
\end{equation}
which can be simplified as
\begin{equation}
\frac{1}{2}\sin\left(2\theta\right)-\theta<0.
\label{eq:dt3}
\end{equation}
Figure \ref{fig:1} presents a plot of the left hand member of \eqref{eq:dt3}, which clearly shows that the inequality is true in the interval of interest. This proves part (b): $d\tau/d\omega<0$.

Finally, statement (c) is a direct conclusion of (a) and (b). Since $\tau$ is a monotonic function of $\omega$ and $\omega$ a monotonic function of $\lambda$, using the chain rule we observe that
\begin{equation}
\frac{d\tau}{d\lambda}=\frac{d\tau}{d\omega}\frac{d\omega}{d\lambda}<0,
\end{equation}
which implies that the largest value of $\lambda$ always invites the minimum $\tau$, and therefore the most restrictive stability boundary. This completes the proof of lemma \ref{le:me2}.
$\qed$

\begin{figure}
\centering
\includegraphics[scale=0.7]{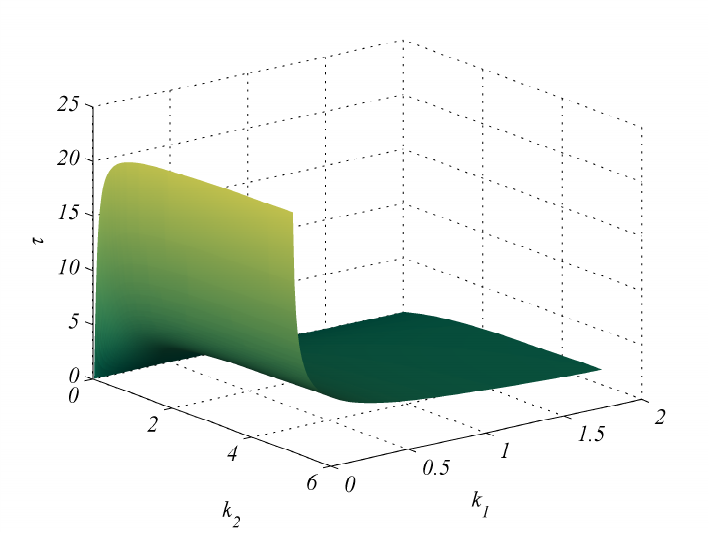}
\caption{Absolute stability boundary in $\left(k_1,\,k_2,\tau\right)$ domain for a group of agents operating under protocol \eqref{eq:noself} with an unknown communication topology. Parametric selections below the surface guarantee stability and convergence to consensus.}
\label{fig:boundarynonself}
\end{figure}

The knowledge of which eigenvalue creates the most restrictive value implies that one does not need to know under which communication topology the agents are operating in order to guarantee stability with respect to the delay. Let us consider a group of agents operating under protocol \eqref{eq:noself}. According to Lemma \ref{le:me1}, the most exigent eigenvalue is the smallest eigenvalue of the weighted adjacency matrix of the communication topology under which the agents are operating. But this eigenvalue is lower bounded: it can never be smaller than $-1$. This implies that the stability boundary defined by \eqref{eq:tau1} for $\lambda=-1$ is the most restrictive for any possible topology. Figure \ref{fig:boundarynonself} shows a plot of this boundary in the tree dimensional domain $\left(k_1,\,k_2,\tau\right)$. Selecting a parametric combination below this surface guarantees that the agents reach consensus \emph{regardless of the topology under which they are operating}, provided that connectivity is present. This is also independent from the number of agents.

In order to perform a similar analysis for protocol \eqref{eq:self}, we need to consider an upperbound for the eigenvalues of the Laplacian of a graph. The first and most conservative of such bounds was presented in \cite{Anderson1985}, and is states that the eigenvalues of the Laplacian of a connected graph are always less than or equal to twice the largest vertex degree of the graph, i.e., $\lambda\leq2\max_i\,\delta_i$. Some other authors have presented less conservative bounds \cite{Li1997,Pan2002}, but they require to have an extra knowledge of the structure of the graph. Without any \emph{a priori} knowledge of the structure of the communication topology, using the bound given by \cite{Anderson1985} is the safest approach.

For a group of $n$ agents the highest possible degree is, of course, $n-1$. One could then consider $\lambda=2\left(n-1\right)$ as the upperbound for the most exigent eigenvalue. However, reference \cite{Anderson1985} states that the upperbound is reached when the graph is bipartite\footnote{A graph is bipartite if its vertices can be separated in two sets such that no vertex is adjacent to another member of its set \cite{Biggs1993}.} and regular\footnote{A graph is regular if all vertices have the same index, i.e., $\delta_i=\delta_j$ for every $i,\,j\in[1,\,n]$ \cite{Biggs1993}}. The highest possible degree for a bipartite graph with $n$ vertices is $n/2$, and therefore the highest Laplacian eigenvalue is given by $\lambda=n$. This is indeed supported by the results presented by \cite{Li1997,Pan2002}.

\begin{figure}
\centering
\subfigure{\includegraphics[scale=0.7]{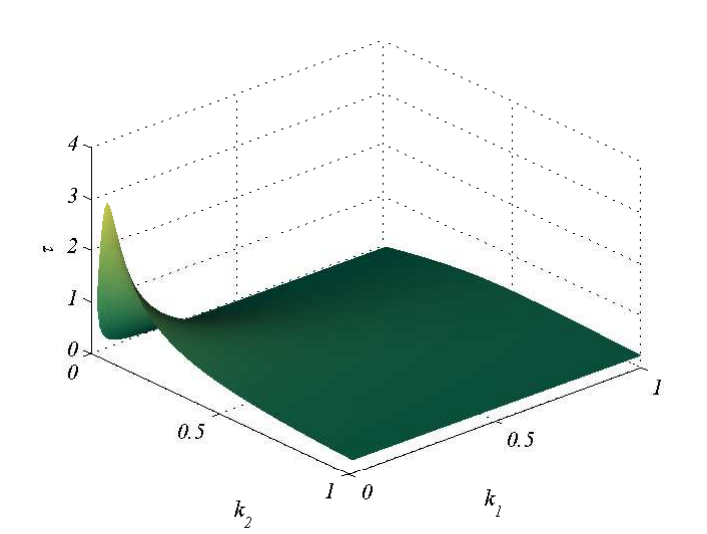}\label{fig:b1self}}
\subfigure{\includegraphics[scale=0.7]{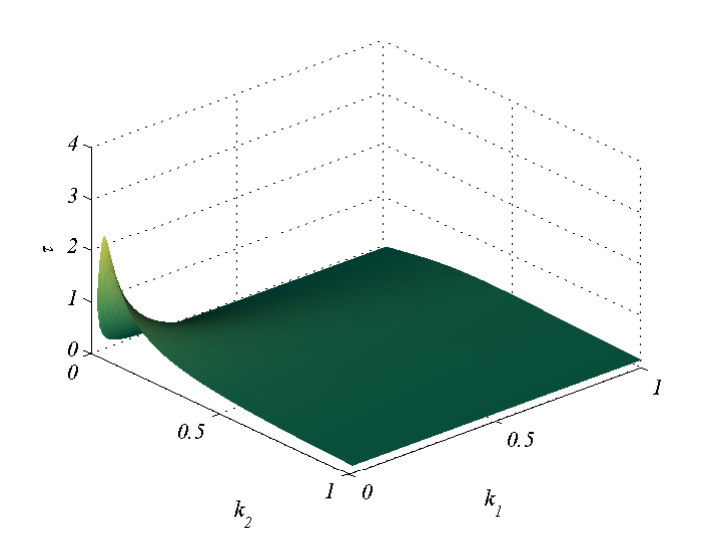}\label{fig:b2self}}
\caption{Absolute stability boundary in $\left(k_1,\,k_2,\tau\right)$ domain for groups of 6 (a) and 10 (b) agents  operating under protocol \eqref{eq:self} with an unknown communication topology. Parametric selections below the surfaces guarantee stability and convergence to consensus.}
\label{fig:boundaryself}
\end{figure}

Since the largest and most exigent Laplacian eigenvalue is dependent on the number of agents, so is the delay margin for a team operating under protocol \eqref{eq:self}. Figure~\ref{fig:b1self} shows the absolute stability boundary for a case with 6 agents, whereas Fig.~\ref{fig:b2self} shows it for a group of 10 agents. Although the difference is not too pronounced, it is possible to see that the boundary is more restricted for larger groups. From an intuitive point of view, this makes sense, because with more communication channels the negative effect of the delay is increased. 

The main driver of this increased effect of the delay is the fact that as the degree of the agents grow, so are the Gershgorin circles for the Laplacian matrix growing. With increasing agents and more connections, the most exigent eigenvalue for this protocol grows and the delay margin decreases. To avoid this, the matrix could be scaled by $n$. Protocol \eqref{eq:noself}, for example, includes this scaling in its definition, when the position and velocity values of the informers are averaged, and therefore its most exigent eigenvalue is well bounded.
  
\section{Switching Topologies}\label{sec:switch}
Section \ref{sec:mee} shows that we can find a parametric combination such that protocols \eqref{eq:noself} and \eqref{eq:self} guarantee consensus regardless of the topology under which the agents are communicating, provided that it is connected. In this section, we focus our attention to the case in which the topologies are switching, and check if by stabilizing all the possible communication topologies the stability of the switching system is guaranteed.

Let us now assume that the communication topology is not fixed, but it switches among a set of $m$ different topologies, described by the weighted adjacency matrices $\vect{C}_i$ and the Laplacians $\vect{L}_i$, with $i=1,\,2,\ldots,\,m$. Each one of these $m$ topologies is assumed to be connected. The switching instants will be denoted as $0<t_{s_{1}}<t_{s_{2}}<\ldots$. There is no assumption on the order of switching.

We revisit now the state space formulation \eqref{eq:genss}. In order to obtain a factorized characteristic equation, we used a state transformation defined by the diagonalization of matrix $\vect{M}$. If $\vect{T}^{-1}\vect{M}\vect{T}=\boldsymbol{\Lambda}$, where $\boldsymbol{\Lambda}$ is a diagonal matrix, the state transformation is defined as $\boldsymbol{\xi}(t)=\left(\vect{T}^{-1}\otimes\vect{I}_2\right)\vect{x}(t)$. The new state vector $\boldsymbol{\xi}(t)$ can be seen as the concatenation of $n$ state vectors $\left[\xi_i(t)\ \dot{\xi}_i(t)\right]^T$, $i=1,\,2,\,\ldots,\,n$, each one corresponding to a subsystem with dynamics
\begin{equation}
\left[\begin{array}{c}\dot{\xi}_i(t)\\\ddot{\xi}_i(t)\end{array}\right]=
\vect{F}_1\left[\begin{array}{c}\xi_i(t)\\\dot{\xi}_i(t)\end{array}\right]+\lambda_i\vect{F}_2\left[\begin{array}{c}\xi_i(t-\tau)\\\dot{\xi}_i(t-\tau)\end{array}\right].
\label{eq:subsys2}
\end{equation}

The characteristic equation of each subsystem of the form \eqref{eq:subsys2} corresponds to one factor $\factor$ in the global characteristic equation \eqref{eq:ce}. For $i=1$ we have the system that creates the centroid factor, and it is therefore called the \emph{centroid subsystem}. We refer to other subsystems in \eqref{eq:subsys} (for $i=2,\,2,\,\ldots,\,n$) as the \emph{disagreement subsystems}. The \emph{disagreement vector} is now defined as a vector containing the state of all the disagreement subsystems
\begin{equation}
\boldsymbol{\xi}_d(t)=\left[\xi_2(t)\,\dot{\xi}_2(t)\,\xi_3(t)\,\dot{\xi}_3(t)\,\cdots\,\xi_n(t)\,\dot{\xi}_n(t)\right].
\label{eq:disagreementv}
\end{equation}
The transformed state vector can be seen as the vector sum of the state of the centroid subsystem and the norm of the disagreement vector, and therefore
\begin{equation}
\left\|\boldsymbol{\xi}(t)\right\|^2=\left\|\boldsymbol{\xi}_1(t)\right\|^2+\left\|\boldsymbol{\xi}_d(t)\right\|^2.
\label{eq:normdecomp}
\end{equation}

For a group of agents operating under a fixed communication topology, a parametric selection within the boundary defined by the most exigent eigenvalue guarantees that all the disagreement subsystems are stable, implying that the norms of their states decrease exponentially with time. The disagreement vector follows this behavior, its norm decreases until consensus is reached. This is formally stated in the following Lemma.

\begin{lemma}\label{le:lyap}
For a set of agents operating under a fixed communication topology, the norm of the disagreement vector is always decaying when the parametric selection is such that the stability of the disagreement subsystems is guaranteed.
\end{lemma}
\emph{Proof:} Let the state of each disagreement subsystem be $\boldsymbol{\xi}_i(t)=[\xi_i(t)\ \dot{\xi}_i(t)]^{T}$, $i=2,\,3,\ldots\,n$. Given that every disagreement subsystem is stable, the Lyapunov-Krasowskii approach \cite{Friedman2014} allows the construction of a Lyapunov function of the form $V_i=\boldsymbol{\xi}_i(t)^{T}\boldsymbol{\xi}_i(t)=\left\|\boldsymbol{\xi}_i(t)\right\|^2$ for which $\dot{V}_i(t)<0$. Since the derivative of the norm of the individual disagreement subsystems is negative definite, these norms are decreasing monotonically. The norm of the full disagreement vector 
\begin{equation}
\left\|\boldsymbol{\xi}_d(t)\right\|=\sqrt{\sum_{i=2}^n{\left\|\boldsymbol{\xi}_i(t)\right\|}}
\end{equation}
is therefore also monotonically decreasing. $\qed$

Lemma \ref{le:lyap} guarantees that for $t_{s_{i}}<t<t_{s_{i+1}}$, \emph{i.e.}, between two switching instants, the norm of the disagreement vector decays and the agents tend to converge into consensus. The end results depends then from what happens right at the switching instants. The following subsections describe what happens in the case of each protocol.

\subsection{Protocol \eqref{eq:noself}}
For a group of agents operating under protocol \eqref{eq:noself}, the group decision value, presented in \eqref{eq:weightedc}, is topology dependent because the degrees of the agents are in general different for each topology. If the topology switches the agents try to arrive to a different decision value, \emph{i.e.}, a jump in the norm of the centroid vector $\left[\xi_1(t)\,\dot{\xi}_1(t)\right]^T$ occurs. This also causes a jump in the norm of the disagreement vector $\boldsymbol{\xi}_d(t)$ at the switching instant, due to \eqref{eq:normdecomp}. Because of the unknown nature of the transformation matrix, there is no way to guarantee that the jumps are bounded. This behavior may induce instability, even when the system is switching among stable topologies.
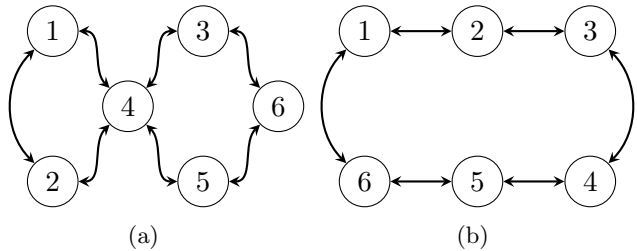
\begin{figure}
\centering
\subfigure[$$]{
\begin{tikzpicture}
\node (n1) [draw,circle] at (-1.5,1) {1};
\node (n2) [draw,circle] at (-1.5,-1){2};
\node (n3) [draw,circle] at (0.5,1)  {3};
\node (n4) [draw,circle] at (-0.5,0) {4};
\node (n5) [draw,circle] at (0.5,-1) {5};
\node (n6) [draw,circle] at (1.5,0)  {6};
\draw [thick,stealth-stealth] (n1.south west) [out=225,in=135] to (n2.north west);
\draw [thick,stealth-stealth] (n1.east) [out=0,in=135] to (n4.north west);
\draw [thick,stealth-stealth] (n2.east) [out=0,in=225] to (n4.south west);
\draw [thick,stealth-stealth] (n3.west) [out=180,in=45] to (n4.north east);
\draw [thick,stealth-stealth] (n5.west) [out=180,in=315] to (n4.south east);
\draw [thick,stealth-stealth] (n3.east) [out=0,in=135] to (n6.north west);
\draw [thick,stealth-stealth] (n5.east) [out=0,in=225] to (n6.south west);
\end{tikzpicture}\label{fig:topol1}}
\subfigure[$$]{\begin{tikzpicture}
\node (n1) [draw,circle] at (-1.5,1) {1};
\node (n2) [draw,circle] at (0,1)    {2};
\node (n3) [draw,circle] at (1.5,1)  {3};
\node (n4) [draw,circle] at (1.5,-1) {4};
\node (n5) [draw,circle] at (0,-1)   {5};
\node (n6) [draw,circle] at (-1.5,-1){6};
\draw [thick,stealth-stealth] (n1.east) to (n2.west);
\draw [thick,stealth-stealth] (n2.east) to (n3.west);
\draw [thick,stealth-stealth] (n6.east) to (n5.west);
\draw [thick,stealth-stealth] (n5.east) to (n4.west);
\draw [thick,stealth-stealth] (n1.south west) [out=225,in=135] to (n6.north west);
\draw [thick,stealth-stealth] (n3.south east) [out=315,in=45] to (n4.north east);
\end{tikzpicture}\label{fig:topol2}}
\caption{Two different connected topologies for six agents.}
\label{fig:topol}
\end{figure}

To better illustrate this idea, consider a group of six agents operating under protocol \eqref{eq:noself} and switching among the two different communication topologies depicted in Fig.~\ref{fig:topol}. With a parametric setting of $\left(k_1,\,k_2,\,\tau\right)=\left(5,\,0.2,\,0.06\right)$ we guarantee stability for each individual topology. We defined the switching pattern by means of a periodic signal with period $T=1.4$ s and a duty cycle $\alpha$, such that during $\alpha \%$ of the period the agents communicate according to the topology depicted in Fig.~\ref{fig:topol1}, and according to the topology in Fig.\ref{fig:topol2} during the rest of it.

\begin{figure}
\centering
\subfigure[]{\includegraphics[scale=0.8]{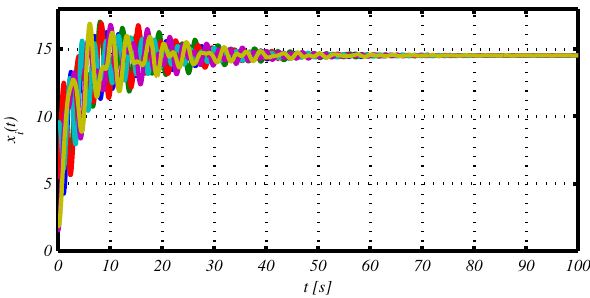}\label{fig:stablea}}
\subfigure[]{\includegraphics[scale=0.8]{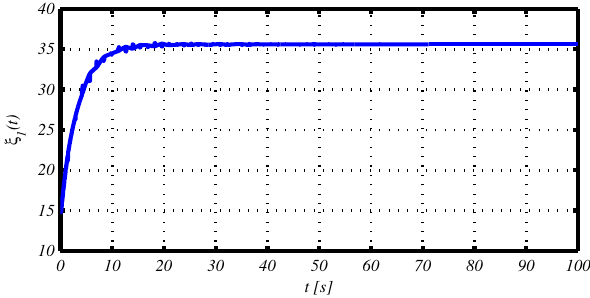}\label{fig:stableb}}
\subfigure[]{\includegraphics[scale=0.8]{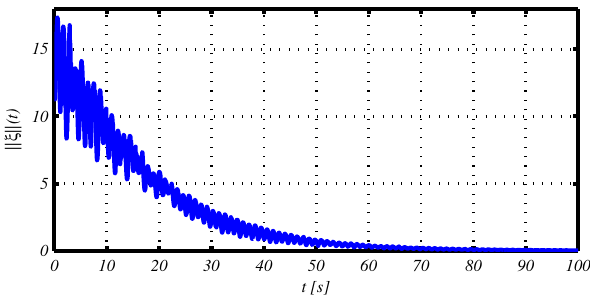}\label{fig:stablec}}
\caption{Behavior of six agents operating under protocol \eqref{eq:noself} and switching stably among the two communication topologies of Fig.\ref{fig:topol}. a) Positions of the agents. b) position of the weighted centroid. c) Norm of the disagreement vector.}
\label{fig:stable}
\end{figure}

Figure \ref{fig:stable} depicts the time evolution of the positions of the agents (Fig.~\ref{fig:stablea}), the topology dependent weighted centroid (Fig.~\ref{fig:stableb}), and the norm of the disagreement vector for a case in which $\alpha=10$. It is clear that, despite the jumps on the value of the weighted centroid, the agents reach consensus.

\begin{figure}
\centering
\subfigure[]{\includegraphics[scale=0.8]{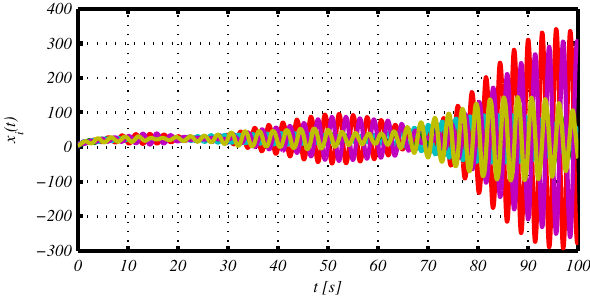}\label{fig:unstablea}}
\subfigure[]{\includegraphics[scale=0.8]{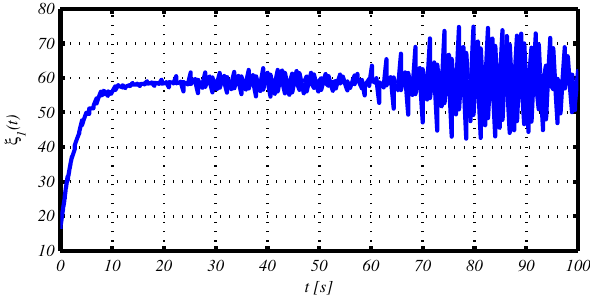}\label{fig:unstableb}}
\subfigure[]{\includegraphics[scale=0.8]{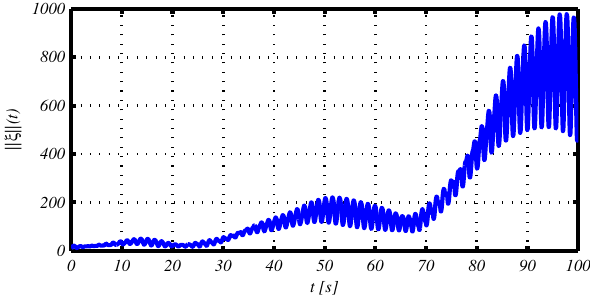}\label{fig:unstablec}}
\caption{Behavior of six agents operating under protocol \eqref{eq:noself} and switching unstably among the two communication topologies of Fig.\ref{fig:topol}. a) Positions of the agents. b) position of the weighted centroid. c) Norm of the disagreement vector.}
\label{fig:unstable}
\end{figure}

By changing only the duty cycle of the switching signal to $\alpha=60$ completely different results are obtained. Figure~\ref{fig:unstable} shows how the positions of the agents, the weighted centroid, and the norm of the disagreement vector evolve under this conditions.

This example shows that when a multi-agent system operates under protocol \eqref{eq:noself} and switching topologies, it is not possible to guarantee stability in the general case. The conditions under which the switching system is stable or unstable remain as an open question for further study.

\subsection{Protocol \eqref{eq:self}}
For protocol \eqref{eq:self} the following lemma proves the most important point.
\begin{lemma}\label{le:cont}
For a group of agents interacting under protocol \eqref{eq:self} the norm of the disagreement vector is continuous at the switching instant. 
\end{lemma}
\emph{Proof:} Since the Laplacian matrix $\vect{L}$ is symmetric, the diagonalization $\vect{T}^{-1}\vect{M}\vect{T}=\boldsymbol{\Lambda}$ is performed using an orthogonal matrix $\vect{T}$, for which $\vect{T}^{-1}=\vect{T}^T$. When an orthogonal matrix multiplies a vector, the norm of the vector remains unchanged. Consider then a switch between two different topologies, described by Laplacians $\vect{L}_1$ and $\vect{L}_2$ which are diagonalized by orthogonal matrices $\vect{T}_1$ and $\vect{T}_2$. The transformed vectors before and after the change in topology are given by
\begin{subequations}\label{eq:swself}
\begin{align}
\boldsymbol{\xi}^{(1)}(t)&=\left(\vect{T}_1^{T}\otimes\vect{I}_2\right)\vect{x}(t),\label{eq:swself1}\\
\boldsymbol{\xi}^{(2)}(t)&=\left(\vect{T}_2^{T}\otimes\vect{I}_2\right)\vect{x}(t).\label{eq:swself2}
\end{align}
\end{subequations}
Given that the agents are not changing their positions or velocities when the communication topology changes, $\left\|\vect{x}(t)\right\|$ is the same in both equations \eqref{eq:swself}. Now, considering that $\vect{T}^T\otimes\vect{I}_2$ is an orthogonal matrix, we have that 
\begin{equation}
\left\|\boldsymbol{\xi}^{(1)}(t)\right\|=\left\|\vect{x}(t)\right\|=\left\|\boldsymbol{\xi}^{(2)}(t)\right\|.
\label{eq:switch1}
\end{equation}
As stated in \eqref{eq:c}, for agents operating under protocol \eqref{eq:self} the state of the centroid subsystem is always equal to the average of the positions of the agents. Since the positions are not changing, their average does not change either,  therefore
\begin{equation}
\left\|\boldsymbol{\xi}^{(1)}_1(t)\right\|=\left\|\boldsymbol{\xi}^{(2)}_1(t)\right\|.
\label{eq:switch2}
\end{equation}
Considering \eqref{eq:normdecomp}, \eqref{eq:switch1} and \eqref{eq:switch2} imply that the disagreement vector $\boldsymbol{\xi}_d$ is continuous at the switching instants. $\qed$

We can now state the final result regarding the stability of a group of agents operating under protocol \eqref{eq:self} and switching topologies.
\begin{lemma}\label{le:final}
If the parametric selection is such that every individual topology leads to consensus, a group of agents operating under protocol \eqref{eq:self} and individual topologies always reaches consensus, regardless of the switching scheme.
\end{lemma}
\emph{Proof:} in order to reach consensus, the norm of the disagreement vector should tend to zero. Lemma \ref{le:lyap} show that this norm decays between switching instants. Lemma \ref{le:cont} shows that it does not change at switching instants. we can conclude that the norm of the disagreement vector is always decreasing, and therefore the agents reach consensus.
$\qed$

\emph{Remark:} This result is equivalent to that found by Olfati-Saber and Murray in \cite{Olfati-Saber2004} for the case of protocols with first order agents. We presented here a non-trivial extension to the case of consensus for second order systems.

The implication of Lemma \ref{le:final} is that if a multi-agent system operating under protocol \eqref{eq:self} switches among stable communication topologies, it reaches consensus regardless of the switching scheme. In this situation, the most exigent eigenvalue is very useful. It allows the designer to select a proper set of parameters which guarantees consensus without knowing the topologies or the switching scheme.

\begin{figure}
\centering
\subfigure[]{\includegraphics[scale=0.8]{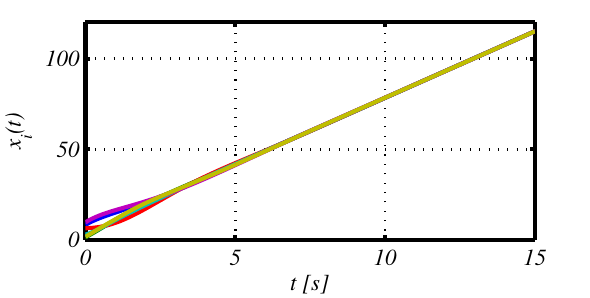}\label{fig:selfa}}
\subfigure[]{\includegraphics[scale=0.8]{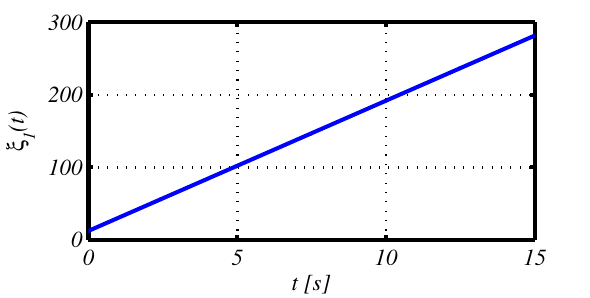}\label{fig:selfb}}
\subfigure[]{\includegraphics[scale=0.8]{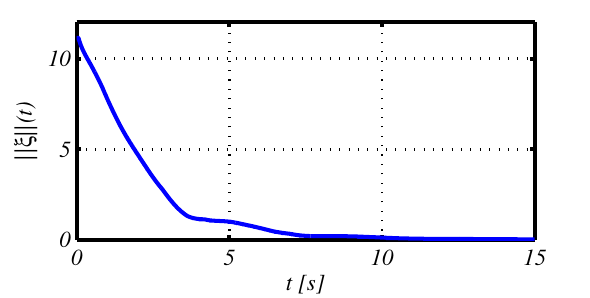}\label{fig:selfc}}
\caption{Behavior of six agents operating under protocol \eqref{eq:self} and switching stably among the two communication topologies of Fig.\ref{fig:topol}. a) Positions of the agents. b) position of the weighted centroid. c) Norm of the disagreement vector.}
\label{fig:self}
\end{figure}

An example of this is presented in Fig.~\ref{fig:self}. It shows the results of a simulation in which six agents operate under protocol \eqref{eq:self} and switch among the two topologies of Fig.~\ref{fig:topol}. The switching is again periodic, and the period and duty cycle are the same as in the unstable example of Fig.~\ref{fig:unstable}. The parameters are selected below the corresponding surface in Fig.~\ref{fig:b1self}. Notice that Fig.~\ref{fig:selfc} shows a smooth trace, confirming the results of Lemma \ref{le:me2}.

\section{Conclusions}\label{sec:conc}
This paper studies two consensus protocols for groups of agents driven by second order dynamics and affected by communication time delays. The agents are assumed to be operating under an bidirectional scheme, such that the communication topology is described by an undirected graph. The relation between the eigenvalues of a matrix associated to such graph and the maximum delay value for which the agents are able to reach consensus is studied in detail for both cases.

It is shown that for one of the protocols there is a possibility of defining a stability boundary that guarantees convergence to consensus even without knowing the number of agents. For the other protocol, the knowledge of the number of agents suffices to define a stability boundary which guarantees consensus for any connected topology under which the agents could be operating. The eigenvalue that defines the most restrictive boundary in each case is called \emph{the most exigent eigenvalue}.

The switching topologies case is also taken under consideration. It is shown that while for one of the protocols the stability can be guaranteed regardless of the switching scheme, for the other one this is not the case and some switching schemes may lead to instability.

Further questions on this topic are related to another concept introduced in \cite{IJC1}: \emph{the most critical eigenvalue}. The question to be answered in this case is: what is the parametric combination $\left(k_1,\,k_2,\,\tau\right)$ which guarantees the fastest convergence to consensus?

The codes used to create the examples of the paper are publicly available in \texttt{https://bit.ly/37pF6Eo}, or can be requested via email to the authors.
%-----------------------------------------------------------------------------------------------------------------
\bibliographystyle{elsarticle-num} 
\bibliography{swarms}
\end{document}